\documentclass[a4,11pt]{article}

\usepackage{amsmath}

\usepackage[english]{babel} 

\usepackage{t1enc}
\usepackage{latexsym}

\DeclareSymbolFont{AMSb}{U}{msb}{m}{n}
\DeclareSymbolFontAlphabet{\mathbb}{AMSb}

\begin{document}  

\bibliographystyle{plain}

\newcommand{\nc}{\newcommand}
\nc{\nt}{\newtheorem}
\nt{defn}{Definition}
\nt{lem}{Lemma}
\nt{pr}{Proposition}
\nt{theorem}{Theorem}
\nt{cor}{Corollary}
\nt{ex}{Example}
\nt{ass}{Assumption}
\nt{step}{Step}
\nt{case}{Case}
\nt{subcase}{Subcase}
\nt{note}{Note}
\nc{\bd}{\begin{defn}} \nc{\ed}{\end{defn}}
\nc{\blem}{\begin{lem}} \nc{\elem}{\end{lem}}
\nc{\bpr}{\begin{pr}} \nc{\epr}{\end{pr}}
\nc{\bth}{\begin{theorem}} \nc{\eth}{\end{theorem}}
\nc{\bcor}{\begin{cor}} \nc{\ecor}{\end{cor}}
\nc{\bex}{\begin{ex}}  \nc{\eex}{\end{ex}}
\nc{\bass}{\begin{ass}}  \nc{\eass}{\end{ass}}
\nc{\bstep}{\begin{step}}  \nc{\estep}{\end{step}}
\nc{\bcase}{\begin{case}}  \nc{\ecase}{\end{case}}
\nc{\bsubcase}{\begin{subcase}}  \nc{\esubcase}{\end{subcase}}
\nc{\bnote}{\begin{note}}  \nc{\enote}{\end{note}}
\nc{\prf}{{\bf Proof.} }
\nc{\eop}{\hfill $\Box$ \\ \\}
\nc{\argmin}{\mathrm{argmin}}
\nc{\argmax}{\mathrm{argmax}}
\nc{\sgn}{\mathrm{sgn}}
\nc{\Var}{\mathrm{Var}}
\nc{\Cov}{\mathrm{Cov}}
\nc{\bak}{\!\!\!\!\!}
\nc{\IBD}{\mathrm{IBD}}
\nc{\supp}{\mathrm{supp}}
\nc{\dom}{\mathrm{dom}}
\nc{\R}{{\mathbb R}}
\nc{\peq}{\preceq}
\nc{\wt}{\widetilde}

\newcommand{\notex}[1]
{$^{(!)}$\marginpar[{\hfill\tiny{\sf{#1}}}]{\tiny{\sf{(!) #1}}}}



\newcommand\approxin
  {\raisebox{-1ex}{ $\stackrel{\textstyle\in}{\scriptstyle\sim}$ }}

\title{Order restricted estimation of the parameter functions in an additive hazard model }
\author{Dragi Anevski\footnote{dragi@maths.lth.se} , Center for Mathematical Sciences, Lund University \\ and \\
ElBatoul Manel Merai\footnote{manel-elbatoul.merai@doc.umc.edu.dz}, Department of Mathematics, \\ Constantine 1 Brothers Mentouri University}\date{\today}
\date{}
\maketitle

\begin{abstract}
In this paper we propose estimators of the parameter functions in an Aalen additive hasard regression model. The estimators are the individual and component wise $l^2$ projections of the naive estimators resulting from the ordinary least squares estimator in the Aalen additive hazard model on the space of monotone functions. The explicit restrictions imposed on the parameter function by our method have a clear clinical interpretation as they can be used to model early onset or late onset covariate effects on complex diseases. We provide point wise limit distribution results for the resulting estimators, that exhibit $n^{1/3}$ rate of convergence and the Chernoff distribution as the limit distribution. We also provide an asymptotic minimax lower bound showing that our suggested estimators are rate optimal, as well as provide an asymptotic minimax lower bound for the more general problem of estimating a monotone hasard function.

 \end{abstract}

\section{Introduction}

In this paper we suggest estimators for the parameter functions in an Aalen additive hasard model, in a survival analysis setting, under the assumption of right-censored data and independent censoring. 

Assuming that the interesting time to event $T$ is a continuous positive random variable, with hasard cumulative distribution function $F$ and hasard function $h(t)=F'(t)/(1-F(t))$, a possible model for $h$, incorporating covariates, is the Aalen additive hasard model 
\begin{eqnarray}\label{eq:Aalen.model}
    h(t)&=&\beta_0(t)+\beta_1(t) z_1 +\ldots+\beta_p(t) z_p,
\end{eqnarray}
where one supposes that  $\beta_0,\ldots,\beta_p$ are (unknown) functions, and $z_1,\ldots,z_p$ are given covariates. If the the parameter vector of function $\beta=(\beta_0,\ldots,\beta_{p+1})$ is completely unspecified, the common approach to estimating the vector of functions $\beta$ is to first realise that it is not possible to provide a nonparametric estimator of it directly. Instead one estimates the vector of integrated functions $B=(B_0,\ldots,B_{p})$ where $B_k(t)=\int _{0}^t \beta_k(u)\,du$, for each $k=0,\ldots,p$, cf. \cite{andersen.borgan.gill.keiding.1993}.

The disadvantage with providing estimators of $B$ instead of $\beta$ is that $B_k(t)$ gives the total effect of covariate $z_k$ summed (i.e.\,integrated) for all times $u\in[0,t]$, whereas $\beta_k(t)$ gives the effect of the covariate $z_k$ {\em at the time $t$}. Clearly $\beta_k(t)$ is more informative and more interesting e.g.\,for the clinical doctor that is interested in describing the effect of covariate values (e.g.\,LDL cholesterol) on the conditional probability of experiencing the interesting event (e.g.\,heart attack) at time $t$, conditional on not having experienced it before time $t$, by the interpretation of the hasard as
\begin{eqnarray*}
    h(t)dt&=&P(T \leq t+dt | T>t).
\end{eqnarray*}

One possibility for obtaining an estimator for $\beta$ from the estimator of $B$ is to use kernel estimators. However, kernel estimators are somewhat ad-doc, and in particular necessitates the choice of a bandwidth. Furthermore, kernel estimators can be seen as a way to provide a {\em smooth estimator} of the unknown parameter, and using such an estimator assumes, implicitly, that the parameter function satisfies a smoothness assumption, e.g.\,that $\beta_k$ is a $C^2$-function. It is not clear what a smoothness assumption on $\beta_k$, such as that it should be a $C^2$-function, corresponds to in a clinical setting.

In this paper we suggest an approach, based on optimising a criterion function under an explicit restriction imposed on the parameter functions, which  solves the above two drawbacks of the kernel estimator. In fact, our suggested estimator does not need an explicit bandwidth choice, since the bandwidth will be {\em data-adaptive}. Furthermore, the model assumption imposed on the parameter function has a clear, and feasible, {\em clinical interpretation}.
 
We propose in this paper an estimator subject to an {\em order restriction} for the $\beta_k$ functions. We define the order restricted nonparametric estimator of $\beta_k$, derive an algorithm for calculating it, and derive limit distributions for the suggested estimator. The algorithm is a version of the well-known minimum lower set algorithm and the limit distribution is the Chernoff distribution, which are commonly featured in order-restricted nonparametric inference.

The suggested model assumption for the $\beta_k$ functions that they are {\em monotone decreasing}, means that the effect of covariate $z_k$ on the hasard function is decreasing with time, or age. This is a feasible assumption for e.g.\,disease models with {\em early onset covariate effects}, such as for Type I diabetes and gene values as covariates. We note that disease models with late onset covariate effects can be modelled assuming that the $\beta_k$ functions are {\em monotone increasing}. Late onset covariate effects may be feasible as disease models with covariates as environmental factors, such as for Type II diabetes and covariates as an environmental or phenotypic measurement. Order-restricted nonparametric estimation of $\beta_k$ under the assumption of monotone increasing is qualitatively similar to the case we treat, and we will state the corresponding results for that case in the sequel. 

A further main advantage of our suggested estimator is that it is data-adaptive in that it uses the implicit bandwidth, imposed by the data. 

There are some previous results on problems that are similar to this inference problem, \cite{huang.2017} uses a method to monotonize the basic estimator in the Aalen model that is different from ours and shows that his estimator is asymptotically equivalent to the standard estimator, while \cite{chung.ivanova.fine.2024} uses a slightly different additive hasard model, which does not seem to include the Aalen model, proposes an order restricted least square estimator and treat mainly computational issues. Order restricted inference of the hasard function $h$ itself under monotonicity assumption, with an emphasis on testing, in a general setting is given in \cite{hall.vankeiligom.2005} and in a survival analysis setting in \cite{bannerjee.wellner.2001}.

Our estimator is to our knowledge new, and our limit distribution results are to our knowledge as well novel.

The paper is organised as follows. In Section \ref{sec:model-setting} we introduce the probabilistic model for the data, as well as the inference problem that we will treat. The model gives rise to a system of stochastic differential equation, and we review the common and well known least squares solution estimator ${\hat B}$ of the integrals $B=(B_0,\ldots,B_{p})$ of the unknown parameter functions $\beta=(\beta_0,\ldots,\beta_p)$. The least squares solution will serve as a starting estimator for the order restricted estimator that we present next. 

Then, in Section \ref{sec:order-restricted-estimator} we present the component-wise least squares projection of the naive estimator arising from the starting estimator presented in Section  \ref{sec:model-setting}, on the space of decreasing functions. These can be written as the derivative of the function $S(\hat{B}_k)$, where $S$ is the least concave majorant map. 

Next, in Section \ref{sec:limit-distributions}, we derive the main results of this paper, which are the limit distributions of the estimators. We start by writing $\hat{B}$ as a sum of the unknown $B$ and a stochastic process $v_n$. We furthermore rescale and localise the estimator $\hat{B}$ which gives rise to a rescaled deterministic term $g_n$ and a rescaled stochastic term $\tilde{v}_n$. In Theorem \ref{thm:limit-distribution-rescaled-proc} we derive the process limit distribution of the $p+1$ dimensional rescaled process $\tilde{v}_n$ to a Gaussian stochastic process $\tilde{v}$ with a certain covariance structure. In Corollary \ref{corr:limit-ditribution-rescaled} we state the resulting component-wise limit distributions for the individual processes $\tilde{v}_{k,n}$, for every $k=0,\ldots,p$. 

Next, in Lemma \ref{lem:limit-tail-bound} we prove a result on a bound on the tail of the process $\tilde{v}_{k,n}$ that ensures that when applying the least concave majorant map $S$ on the process $g_{k,n}+\tilde{v}_{k,n}$, the tail behaviour of that process will not affect the application of the map $S$ around the origin. In Lemma \ref{lem:limitprocess-bound} we state the analog bound for the tail of the limit process $\tilde{v}_{k}$ that ensures the same thing for the limit process $-s^2+\tilde{v}(s)$, where is fact $-s^2$ is proportional to the uniform limit of $g_{k,n}(s)$.

Then in Theorem \ref{thm:limit-distyribution-cumulative} we state one of the two main results of this paper, namely that the integral $\tilde{B}_k$ of the proposed estimators converges to a limit random variable, as
\begin{eqnarray*}
      n^{2/3} c(t_0)(\tilde{B}_k(t_0)-B_k(t_0)&\stackrel{d}{\to}& S(-s^2+w(s))(0),
  \end{eqnarray*}
with $w(s)$ a two-sided Brownian motion. The constant $c(t_0)$ is specified in Theorem \ref{thm:limit-distyribution-cumulative}.

Next, in Theorem \ref{th:limit-distribution-grenander}, we state the second main result of the paper, namely that the proposed estimator $\tilde{\beta}_k$ converges to a limit random variable, as 
\begin{eqnarray*}
      n^{1/3} c(t_0)(\tilde{\beta}_k(t_0)-\beta_k(t_0))&\stackrel{d}{\to}& S(-s^2+w(s))'(0).
  \end{eqnarray*}
  We note that the rate is $n^{1/3}$, which is common in nonparametric order restricted inference, and that the limit distribution $S(-s^2+w(s))'(0)$ is (proportional to) the Chernoff distribution, which is common in nonparametric order restricted inference.  In Theorem \ref{th:limit-distribution-grenander-2} we state the corresponding limit distribution result for the proposed estimator in the case of a monotone increasing parameter function $\beta_k$.
  
In Section \ref{sec:asymptotic-lower-bound} we provide a asymptotic lower bound for the minimax risk, in Theorem \ref{thm:minimax-bound-beta}, showing that our derived estimator is rate optimal. In this section we also provide an asymptotic lower bound for the problem of estimating a monotone hasard function, in Theorem \ref{thm:lower-bound-0}, which may be of independent interest. 
  
Finally in Section \ref{sec:discussion} we discuss the derived results.

\section{The  survival analysis model setting}\label{sec:model-setting}
Let $T\geq 0$ be a positive continuous random variable with an unknown distribution function $F$. We assume that $T$ models the time to an event. We assume no left-truncation for the data, and that we have the standard right-censoring, i.e. that we observe the minimum of the time $T_i$ and a censoring time $C_i$, together with an indicator for the time being exact, $\delta_i=1\{T_i\leq C_i\}$.  

Introduce the individual counting processes $N_i(t)=1\{t_i\leq t, \delta_i=1\}$ for which one has the stochastic differential equation
\begin{eqnarray}
   dN_i(t)&=& Y_i(t)h(t)dt+ dM_i(t), \label{eq:SDE-base}
\end{eqnarray}
for $i=1,\ldots,n$, where $h(t)$ is the individual hasard function, which exists when the distribution of $T$ is absolutely continuous, and is then given by $h(t)=-{d}/{dt}\log(1-F(t))$, where $Y_i(t)=1\{t_i\geq t\}$ is the left-continuous indicator process for the individual being at risk at time $t-$, and $M_i$ is the individual martingale, i.e. satisfying  $E(dM_i(t)| {\cal F}_t)=0$, and where the sigma algebras $\{{\cal F}_t, t\geq \}$ is a filtration storing the information available at times $t$.

The $\sigma$-algebra generated by the information depends on the amount and type of information about the observed times that is available and the amount of information about the covariates that is available. Starting with the model $(\ref{eq:SDE-base})$ satisfying $E(dM_i(t)| {\cal F}_t)=0$ one needs to establish that $E(dM_i(t)| {\cal G}_t)=0$ where ${\cal G}_t$ is the $\sigma$-algebra generated by the observed and available information at time $t$. The concept of noninformative and independent censoring as well as the innovation theorem are used to establish this link. We do not make this assumptions explicit here, since they are of no relevance to us, and refer to reader to a standard reference such as \cite{andersen.borgan.gill.keiding.1993}. We will in the sequel assume that ${\cal F}_t$ is a filtration, containing the information available at time $t$, and that $E(dM_i(t)| {\cal F}_t)=0$ holds so that $E(dN_i(t)| {\cal F}_t)=Y_i(t)h(t)dt$.

A basic inference problem in survival analysis is assessing the effect of group indicators or continuous covariate measurements on the distribution  of the time to an event. It is then necessary to assume that and to model for the distribution function depending on covariates $z_1,\ldots,z_p$. This can be done using various models. One standard model is the Aalen additive hasard model
\begin{eqnarray}\label{eq:Aalen.model}
    h(t)&=&\beta_0(t)+\beta_1(t) z_1 +\ldots+\beta_p(t) z_p,
\end{eqnarray}
where $\beta_0,\ldots,\beta_p$ are (unknown) functions. Thus $\beta=(\beta_0,\ldots,\beta_p)$ is the unknown parameter vector of functions and $E^{p+1}$  is the parameter space, where $E=\{g:[0,\infty)\to {\mathbb R}, \int g<\infty \}$ is the set of integrable functions on $[0,\infty)$.

The value of $\beta_k(t)$  describes the time-instanteneous effect of covariate $z_k$ on the hasard function $h(t)$.  The standard approach for estimating the $\beta_k$'s is to first acknowledge that they are not possible to estimate directly. Rather one estimates their integrals $B_k(t)=\int_{0}^t \beta_k(u)du$. In fact, one can write $(\ref{eq:SDE-base})$ for the Aalen model as
\begin{eqnarray}
   dN_i(t)&=& Y_i(t)\left\{dB_0(t)+dB_1(t)z_{1i}+\ldots +dB_p(t)z_{pi}\right\}+ dM_i(t), \label{eq:SDE}
\end{eqnarray}
for $i=1,\ldots,n$, where $Y_i(t)$ is the individual at-risk process and $M_i$ is a continuous time martingale, for individual $i=1,\ldots,n$. The $n$ equations $(\ref{eq:SDE})$ can be written on the matrix formulation
\begin{eqnarray*}
   dN(t)&=&Y(t)dB(t)+dM(t),
\end{eqnarray*}
where $B(t)=(B_0(t),\ldots, B_{p}(t))^t$ is the vector of unknown functions, and $Y$ is a $n\times (p+1)$ matrix, with the $i$'th row of $Y$ being 
$Y_i(t)(1, z_{1i},\ldots, z_{pi})$.  

If $J(t)$ is the (predictable) indicator that $Y(t)$ has full rank, and
\begin{eqnarray*}
Y^{-}(t) &=& (Y(t)^{T}Y(t))^{-1} Y(t)^{T}
\end{eqnarray*}
is a (generalised) inverse, then $B$ can be estimated by Aalen's ordinary least squares solution
\begin{eqnarray}\label{eq:aalen-least-squares}
   \hat{B}(t)&=&\int_{0}^t J(u)Y^{-}(u)dN(u).
\end{eqnarray}


The interpretation of a value of the integral $B_k(t)$ is less intuitive than the interpretation of the value of $\beta_k(t)$, since $\beta_k(t)$ is the instantaneous effect, at time $t$, of the covariate $z_k$ of the total hasard $h(t)$ at time $t$. Thus one would really like to get an estimate of $\beta_k$, and this is not possible to do directly. One possibility for estimation of $\beta_k$ itself would be to do kernel smoothing, with the drawback that this is a slightly ad hoc method, with a bandwidth  the user has to specify. Thus it is not automated, or data adaptive.
 
 An alternative for estimation of $\beta_k$, which does not necessitate bandwidth choices, and is data adaptive, is to use estimation under some nonparametric restrictions.  In this paper we suggest to estimate $B$ under the assumption that each $\beta_k$ is a {\em nonincreasing} function, i.e. a function that is (not necessarily strictly) decreasing. This is an order restricted inference problem, and a nonparametric such.

\section{The order restricted estimator}\label{sec:order-restricted-estimator}
We define the order restricted estimators of the $\beta_k$'s as the least squares projection of the increments of the components on the Aalen estimator $\hat{B}$ on the space of monotone functions.  Thus let $\hat{B}_k$ be the $k$'th component in $\hat{B}$, for $k=0,\ldots,p$, and suppose that $\hat{B}_k$ (which is a step function) has $L_k$ incremental steps, at points $t_1,\ldots,t_{L_k}$. Thus $(\Delta_1 \hat{B}_k,\ldots, \Delta_{L_k}\hat{B}_k)$ is the vector of increments, where $ \Delta_{j}\hat{B}_k=\hat{B}_k(t_{j})-\hat{B}_k(t_{j-1})$.  Then we may introduce the naive vector of estimates $\hat{\beta}^{(k)}=(\hat{\beta}_1^{(k)},\ldots,\hat{\beta}_{L_k}^{(k)})$, where $\hat{\beta}_i^{(k)}={ \Delta_{j}\hat{B}_k}/{\Delta_j t}$, where $\Delta_j t= t_j-t_{j-1}$, and  with the convention $t_0=0$.

For any integer number $L$, let ${\cal R}_L=\{\gamma \in {\mathbb R}^L: \gamma_1\geq \ldots \geq \gamma_L\}$ be the set of real vectors that have non-increasing coordinates.  We then define the isotonic regression $\tilde{\beta}^{(k)}$ of $\hat{\beta}^{(k)}$ as
\begin{eqnarray}\label{eq:def-isotonic-regression}
     \tilde{\beta}^{(k)}&=&\argmin_{\gamma \in {\cal R}_{L_k}} \sum_{i=1}^{L_k} ( \hat{\beta}_i^{(k)}- \gamma_i)^2.
\end{eqnarray}
Finally we define the order restricted estimator $\tilde{\beta}_k$ as the constant interpolation of the vector $\tilde{\beta}^{(k)}$
\begin{eqnarray*}
    \tilde{\beta}_k(s)&=&\sum_{t_i\leq s}  \tilde{\beta}^{(k)}_i \Delta_i t +  \tilde{\beta}^{(k)}_{j} (s-t_{j}), 
\end{eqnarray*}
where $j=\sup\{i: t_i\leq s\}$ is index of the largest $t_i\leq s$. Standard theory for isotonic regression shows that the vector $\tilde{\beta}^{(k)}$, and therefore $\tilde{\beta}_k(s)$, exists. Furthermore, a geometric characterisation of the solution $\tilde{\beta}_k(s)$ is given by
\begin{eqnarray}\label{eq:algorithm-isotonic-regression}
       \tilde{\beta}_k(t)&=&\frac{d}{dt}(S(\hat{B}_k(t)))
\end{eqnarray}
where $S$ is the least concave majorant map and $d/dt$ denotes the left hand derivative. We note also that the corresponding cumulative function
\begin{eqnarray}\label{eq:cumulative-isotonic-regression}
        \tilde{B}_k(t)&=&S(\hat{B}_k(t))
\end{eqnarray}
is an order restricted estimator of the cumulative function $B_k$, and that it is {\em concave}.

\section{The limit distribution results for the estimators}\label{sec:limit-distributions}

We first see that we can write that the $p+1$-dimensional vector of estimators $\hat{B}$ as
\begin{eqnarray*}
\hat{B}(t) &= &\int_{0}^{t} J(u) Y^{-}(u) dM(u),\\
\end{eqnarray*}
where $J(t)$ is the (predictable) indicator that $Y(t)$ has full rank, where 
\begin{eqnarray*}
Y^{-}(t) &=& (Y(t)^{T}Y(t))^{-1} Y(t)^{T}
\end{eqnarray*}
is a (generalised) inverse, and where $M(t)$ is an $n$-vector of locally square-integrable martingales.

We center the estimator $\hat{B}$ at $B$, and define the process part of the estimator as
\begin{eqnarray*}
v_n(t)&=&
 \int_{0}^{t} J(u) Y^{-}(u) dM(u)
- B(t).
\end{eqnarray*}
to get
\begin{eqnarray*}
\hat{B}(t) &= & B(t)+ v_n(t),
\end{eqnarray*}
and note that this is a slight adaptation from the partition/centering used in \cite{anevski.hossjer.2006}. In fact, we have written the preliminary estimator $\hat{B}(t)$ as a sum of a deterministic part $B(t)$ and a stochastic part $v_n(t)$. The final order restricted estimator is obtained as a coordinate-wise isotonic regression of the increments of the preliminary estimator $\hat{B}(t)$, as defined in $(\ref{eq:def-isotonic-regression})$. 

Therefore the local rescaling should be defined coordinate-wise, and it is in fact enough to study the coordinate-wise partition
\begin{eqnarray*}
\hat{B}_k(t) &= & B_k(t)+ v_{k,n}(t)
\end{eqnarray*}
and the component-wise rescaling
\begin{eqnarray*}
       \tilde{v}_{k,n}(s)&=&d_{n}^{-2} \big( v_{k,n}(t_{0} + s d_{n}) - v_{k,n}(t_{0}) \big),
 \end{eqnarray*}
and to establish limit properties for the rescaled $\tilde{v}_{k,n}$, for the results that we will develop here. In particular we want to establish local limit distribution results as well as certain truncation properties for the rescaled process $\tilde{v}_{k,n}$. We are however able to establish the limit distributions for the full vector valued rescaled process, and since that result may be of independent interest, we will state this result. The coordinate wise property will then be a corollary of that result.

We will make frequent referencing to the Cramér-Wold device for processes, that $W_n$, a $d$-dimensional stochastic process converging in distribution to a $d$-dimensional Gaussian $W$ process, is equivalent to the weak convergence of the one-dimensional process $\alpha_1 W_{1,n}+ \ldots +\alpha_d W_{d,n}$ to $\alpha_1 W_{1}+ \ldots + \alpha_d W_{d}$, for every choice of of $\alpha_1,\ldots,\alpha_d$. 

In fact we are going to adapt the proof of Theorem VII of \cite{andersen.borgan.gill.keiding.1993}  which states that the $p+1$-dimensional process $v_n$ converges to a Gaussian process, say $v$, with a certain covariance structure, to our settings. This implies, by the Cramér-Wold device and since a Gaussian process is determined by it's expectation and covariance function, that the $k$'th coordinate $v_{k,n}$ will converge to a Gaussian process $v_k$ with a covariance structure determined from the covariance structure of the full process $v$. One could therefeore rescale the $k$'th coordinate process $v_{k,n}$ and establish limit distributions for that process. As already mentioned, we will instead rescale the full process and invoke the Cramér-Wold device subsequentaly.

Thus let us define the full rescaled process part
\begin{eqnarray*}
\tilde{v}_{n}(s) = d_{n}^{-2} \big( v_{n}(t_{0} + s d_{n}) - v_{n}(t_{0}) \big).
\end{eqnarray*}

For the local limit distribution results we will adapt the proof of Theorem VII.4.1 of \cite{andersen.borgan.gill.keiding.1993} to our settings, and we will establish the limit distribution result under the same assumptions as those in Theorem VII.4.1 Thus we define
 for $j,k,l = 0,1,\ldots,p$, the functions
\begin{eqnarray*}
R^{(1)}_{j}(t) &=& \sum_{i=1}^{n} Y_{i}(t) Z_{ij}(t), \\
R^{(2)}_{jk}(t) &=& \sum_{i=1}^{n} Y_{ij}(t) Y_{ik}(t), \\
R^{(3)}_{jkl}(t) &=& \sum_{i=1}^{n} Y_{ij}(t) Y_{ik}(t) Y_{il}(t).
\end{eqnarray*}

\noindent
Let $0<\kappa<\infty$ be arbitrary, so that $[0,s']$ is an arbitrary compact set.
\bass\label{Ass:1}
For all $j,k,l = 0,1,\ldots,p$, there exist continuous functions 
$r^{(1)}_{j}, r^{(2)}_{jk}, r^{(3)}_{jkl}$ such that as $n \to \infty$:
\begin{align*}
\sup_{s \in [0,\kappa]} \left| \frac{1}{n} R^{(1)}_{j}(s) - r^{(1)}_{j}(s) \right| &\xrightarrow{P} 0,\\
\sup_{s \in [0,\kappa]} \left| \frac{1}{n} R^{(2)}_{jk}(s) - r^{(2)}_{jk}(s) \right| &\xrightarrow{P} 0,\\
\sup_{s \in [0,\kappa]} \left| \frac{1}{n} R^{(3)}_{jkl}(s) - r^{(3)}_{jkl}(s) \right| &\xrightarrow{P} 0.
\end{align*}
\eass

\bass\label{Ass:2}
For all $j=0,1,\ldots,p$,
\begin{eqnarray*}
\frac{1}{\sqrt{n}} \sup_{i=1,\dots,n} \sup_{s \in [0,\kappa]} |Y_{ij}(s)| \xrightarrow{P} 0.
\end{eqnarray*}

\eass
\bass\label{Ass:3}
For all $s \in [0,\kappa]$, the matrix $r^{(2)}(s) = \big( r^{(2)}_{jk}(s) \big)$ is nonsingular.
\eass

\bth\label{thm:limit-distribution-rescaled-proc}
Suppose that Assumptions \ref{Ass:1}- \ref{Ass:3} hold. Then
\begin{eqnarray*}
      \tilde{v}_n(s)&\stackrel{d}{\to}& \tilde{v}(s)
\end{eqnarray*}
on $D^{p+1}(-c,c)$, as $n\to\infty$, where $\tilde{v}$ is mean zero Gaussian process with covariance structure 
\begin{eqnarray*}
Cov( \tilde{v}_{j}(s'), \tilde{v}_{k} (s'')&=&\sigma_{j,k}\, \min(s',s''),
\end{eqnarray*}
where 
\begin{eqnarray*}
 \sigma_{j,k}&=&\sum_{g,l,m=0}^{p}
(r^{(2)}(t_{0}))^{-1}_{jl}
(r^{(2)}(t_{0}))^{-1}_{km}
r^{(3)}_{lmg}(t_{0}) \beta_{g}(t_{0}) .
\end{eqnarray*}
\eth
\noindent
\prf
Defining the matrix 
\begin{eqnarray*}
R^{(2)}(t) &:=& \sum_{i=1}^{n} Y_{i}(t) Y_{i}^{T}(t) \\
&=& Y(t)^{T}Y(t),
\end{eqnarray*}
the second statement of Assumption \ref{Ass:1} implies
\begin{eqnarray*}
\sup_{s\in[0,\kappa]}||\frac{1}{n} R^{(2)}(s)-r^{(2)}(s)||&\stackrel{P}{\to}&0
\end{eqnarray*}
where $r^{(2)}$ is defined in Assumption \ref{Ass:3}, and with $|| \cdot ||$ denoting the euclidian (matrix) norm on ${\mathbb R}^{p+1}\times {\mathbb R}^{p+1}$.  By Assumption \ref{Ass:3} the matrix $r^{(2)}$ is invertible, and the inverse $(\frac{1}{n} R^{(2)}(s))^{-1}$ is well defined when $J(s)=1$,  and for those $s$ converges to $(r^{(2)}(s))^{-1}$ by the continuous mapping theorem, since matrix inversion is a continuous map (under the supnorm matrix metric). Thus, since $Y^{-}(t) = (Y(t)^{T}Y(t))^{-1} Y(t)^{T}$, we may partition the process part as 
\begin{eqnarray}
v_{n}(t) &= & \int_{0}^{t} J(u) \Big[
    \big( \tfrac{1}{n} R^{(2)}(u) \big)^{-1} - \big( r^{(2)}(u) \big)^{-1}
\Big] Y^{T}(u) dM(u) \nonumber \\
&& + \int_{0}^{t} J(u) \big( r^{(2)}(u) \big)^{-1} Y^{T}(u) dM(u)
+ \int_{0}^{t} (J(u) - 1) \beta(u) du. \label{eq:process-partition}
\end{eqnarray}

Recall the definition of the rescaled process
\begin{eqnarray*}
\tilde{v}_{n}(s) = d_{n}^{-2} \big( v_{n}(t_{0} + s d_{n}) - v_{n}(t_{0}) \big), 
\quad
\end{eqnarray*}
for $ s \in [-c, c]$, and note that it entails that $\tilde{v}_{n}(s)$ is the sum of the three integrals in $(\ref{eq:process-partition})$ with the integrals going from $t_0$ to $t_0+sd_n$, and all multiplied by $d_n^{-2}$. We may now use a change of variables inside the integrals,
so for $u \in [ t_{0}, t_{0}+s d_{n}] $ and $s$ fixed we let $ s' \in [ 0, s]$ vary and 
and thus $du= d_{n}ds'$ so that we obtain
\begin{eqnarray*}
\tilde{v}_{n}(s) &= &
\frac{d_{n}^{-2}}{\sqrt{n}} 
    \frac{1}{\sqrt{n}} 
    \int_{0}^{s} J(t_{0}+s'd_{n})
    \Big( \big( \tfrac{1}{n} R^{(2)}(t_{0}+s'd_{n}) \big)^{-1}
        - \big( r^{(2)}(t_{0}+s'd_{n}) \big)^{-1} \Big) \cdot\\
 &&\cdot Y^{T}(t_{0}+s'd_{n}) dM(t_{0}+s'd_{n}) \\
&&+ \frac{d_{n}^{-2}}{\sqrt{n}} 
\left[  \frac{1}{\sqrt{n}} 
    \int_{0}^{s} J(t_{0}+s'd_{n})
    (r^{(2)}(t_{0}+s'd_{n}))^{-1}
    Y^{T}(t_{0}+s'd_{n}) dM(t_{0} +s'd_{n}) \right] \\
&&+ \frac{d_{n}^{-2}}{\sqrt{n}} 
\left[ \sqrt{n} \int_{0}^{s} (J(t_{0}+s'd_{n}) - 1) \beta(t_{0}+s'd_{n})d_{n} ds' \right]\\
&=:&\tilde{v}_{n}^{(1)}(s)+\tilde{v}_{n}^{(2)}(s)+\tilde{v}_{n}^{(3)}(s).
\end{eqnarray*}

We will now treat the three terms $\tilde{v}_{n}^{(1)},\tilde{v}_{n}^{(2)},\tilde{v}_{n}^{(3)}$ above separately, and show that the first and last vanish asymptotically, while the second $\tilde{v}_{n}^{(2)}$ gives rise to the asymptotic distribution. 

$(i): \tilde{v}_{n}^{(1)}$ {\em vanishes asymptotically}. 

If we denote write the $j$'th component of the first term $ \tilde{v}_{n}^{(1)}$ as $\tilde{v}_{n,j}^{(1)}$, we get
\begin{eqnarray*}
\tilde{v}_{n,j}^{(1)}(s)&=&\frac{d_{n}^{-2}}{\sqrt{n}} \frac{1}{\sqrt{n}} 
\int_{0}^{s} J(t_{0}+s'd_{n})
\sum_{i=1}^{n} \sum_{l=0}^{p}
\Big[
    \big( \tfrac{1}{n} R^{(2)}(t_{0}+sd_{n}) \big)^{-1}
   \\
   && - \big( r^{(2)}(t_{0}+s'd_{n}) \big)^{-1}
\Big]_{jl} Y_{il}(t_{0}+s'd_{n}) dM_{i}(t_{0}+s'd_{n}),
\end{eqnarray*}
for arbitrary but fixed $j=0,\ldots,p$.
Therefore the predictable variation process becomes
\begin{eqnarray*}
\langle \tilde{v}_{n,j}^{(1)} \rangle (s) 
&= &\frac{d_{n}^{-4}}{n} \frac{1}{n}
\int_{0}^{s} J(t_{0}+s'd_{n})
\sum_{i=1}^{n}
\Bigg\{ \sum_{l=0}^{p} 
    \Big( \big( \tfrac{1}{n} R^{(2)}(t_{0}+s'd_{n}) \big)^{-1}  \\
    && - \big( r^{(2)}(t_{0}+s'd_{n}) \big)^{-1}
    \Big)_{jl} \cdot Y_{il}(t_{0}+s'd_{n})
\Bigg\}^{2} d\langle M_{i} \rangle (t_{0}+s'd_{n})\\
&=&\frac{d_{n}^{-4}}{n} \frac{1}{n}
\int_{0}^{s} J(t_{0}+s'd_{n})
\sum_{i=1}^{n}
\Bigg\{ \sum_{l=0}^{p} 
    \Big( \big( \tfrac{1}{n} R^{(2)}(t_{0}+s'd_{n}) \big)^{-1}  \\
    && - \big( r^{(2)}(t_{0}+s'd_{n}) \big)^{-1}
    \Big)_{jl} \cdot Y_{il}(t_{0}+s'd_{n}) 
\Bigg\}^{2}  d_n \lambda_{i}(t_0) ds',
\end{eqnarray*}
since $d\langle M_{i} \rangle (t_{0}+s'd_{n})=d_n \lambda_{i}(t_0) ds' $, and since $d \langle M_{i},M_j \rangle=0$ for $i\neq j$. From the above we see that the local rescaling rate $d_{n} = n^{-\alpha}$ is determined by the condition
\begin{eqnarray*}
\frac{d_{n}^{-4}}{n} \cdot d_{n} &=& \frac{d_{n}^{-3}}{n} = 1 
\end{eqnarray*}
and thus we must have $d_{n} = n^{-1/3}$.

By Assumption \ref{Ass:1} we have that 
\begin{eqnarray*}
 \sup\limits_{s' \in [0,s]} || \displaystyle\frac{1}{n} R^{(2)}(t_{0}+s'd_{n}) -r^{(2)}(t_{0}+s'd_n) || &\stackrel{P}{\rightarrow} &0
 \end{eqnarray*}
and by Assumption \ref{Ass:3}, $ r^{(2)}(t_{0})$ is nonsingular. Then, since on the set of points $t$ where $J(t)=1$ we have that $R^{(2)}(t)^{-1}$ exists,    and since the matrix inverse map is a continuous map (under the supnorm metric), then by the continuous mapping theorem, for $j$ fixed and for every $l=0,\ldots,p$,
\begin{eqnarray*}
   \sup_{s'\in [0,s]} | J(t_{0}+sd_{n}) 
\Big(
    (\tfrac{1}{n} R^{(2)}(t_{0}+s'd_{n}))^{-1}
    - (r^{(2)}(t_{0}+s'd_{n}))^{-1}
\Big)_{jl} | &\stackrel{P}{\rightarrow} &0.
 \end{eqnarray*}
Therefore, there are random variables $C_{0}^{(n)},\ldots,C_{p}^{(n)} $ such that
\begin{eqnarray*}
&&\sup_{s'\in [0,s]} | J(t_{0}+s'd_{n}) \sum_{l=0}^{p} 
\Big(
    (\tfrac{1}{n} R^{(2)}(t_{0}+s'd_{n}))^{-1}
    - (r^{(2)}(t_{0}+s'd_{n}))^{-1}
\Big)_{jl}
Y_{il}(t_{0}+s'd_{n})| \\
&\leq&\sum_{l=0}^p C_{l}^{(n)} Y_{il}(t_{0}+s'd_{n}),
\end{eqnarray*}
and such that $C_{l}^{(n)}=o_P(1)$, for all $l=0,\ldots,p$.

Thus, for $d_{n} = n^{-1/3}$,
\begin{eqnarray*}
\langle \tilde{v}_{n,j}^{(1)} \rangle (s) 
&\leq &
\int_{0}^{s} J(t_{0}+s'd_{n})
\sum_{l,k=0}^{p} C_{l}^{(n)} C_{k}^{{n}} 
\frac{1}{n} \sum_{i=1}^{n} 
Y_{il}(t_{0}+s'd_{n}) Y_{ik}(t_{0}+s'd_{n})
\lambda_{i}(t_{0}+s'd_{n}) ds'.
\end{eqnarray*}
From the second line of Assumption \ref{Ass:1}, we have
\begin{eqnarray*}
P\left(
    \sup_{s' \in [0,s]} 
    \Big|
        \frac{1}{n} \sum_{i=1}^{n}
        Y_{il}(t_{0}+s'd_{n}) Y_{ik}(t_{0}+s'd_{n})
        - r^{(2)}_{kl}(t_{0}+s'd_{n})
    \Big| > \varepsilon
\right) \to 0,
\end{eqnarray*}
and thus
\begin{eqnarray*}
\langle \tilde{v}_{n,j}^{(1)} \rangle (s) \xrightarrow{P} 0.
\end{eqnarray*}
Finally by Lenglart’s inequality, for all $\varepsilon, \eta > 0$,
\begin{eqnarray*}
P\left( 
    \sup_{s' \in [0,s]} |\tilde{v}_{n,j}^{(1)}(s')| > \varepsilon 
\right)
&\leq & 
\frac{\eta}{\varepsilon^{2}} + 
P\big( \langle \tilde{v}_{n,j}^{(1)} \rangle (s) > \eta \big),
\end{eqnarray*}
which shows that that $\tilde{v}_{n,j}^{(1)}$ converges in probability to zero, uniformly on $[0,s]$, as claimed.

$(ii): \tilde{v}_{n}^{(2)}$ {\em gives the limiting behaviour.}

We derive the asymptotic normality of $\tilde{v}_{n}^{(2)}$ by first establishing the limit in probability of the predictable covariation processes of $\tilde{v}_{n}^{(2)}$ and then by checking the Lindeberg condition for $\tilde{v}_{n}^{(2)}$.

Let $\tilde{v}_{n,j}^{(2)} $ denote the $j$'th component of $\tilde{v}_{n}^{(2)}$, for $j=0,\ldots,p$. Thus
\begin{eqnarray}
\tilde{v}_{n,j}^{(2)}(s)& =&
\frac{d_{n}^{-2}}{\sqrt{n}} \frac{1}{\sqrt{n}}
\int_{0}^{s} 
\sum_{i=1}^{n} 
\sum_{l=0}^{p}
(r^{(2)}(t_{0}+s'd_{n}))^{-1}_{jl}
Y_{il}(t_{0}+s'd_{n}) dM_{i}(t_{0}+s'd_{n}). \label{eq:v_n-tilde-j}
\end{eqnarray}

We first establish the asymptotic limit of the quadratic covariation processes. For $d_{n} = n^{-1/3}$, the predictable quadratic covariation between components $j$ and $k$ is
\begin{eqnarray*}
\langle \tilde{v}_{n,j}^{(2)}, \tilde{v}_{n,k}^{(2)} \rangle (s) &=&
 \frac{d_{n}^{-4}}{n} \frac{1}{n} \sum_{i=1}^{n}
\int_{0}^{s} 
\sum_{l,m=0}^{1}
(r^{(2)}(t_{0}+s'd_{n}))^{-1}_{jl}
(r^{(2)}(t_{0}+s'd_{n}))^{-1}_{km} \\
&&
\cdot Y_{il}(t_{0}+sd_{n}) Y_{im}(t_{0}+s'd_{n})d \left\langle M_{i} \right\rangle (t_{0}+s'd_{n}) \\
&=&\frac{1}{n} \sum_{i=1}^{n}
\int_{0}^{s} 
\sum_{l,m=0}^{1}
(r^{(2)}(t_{0}+s'd_{n}))^{-1}_{jl}
(r^{(2)}(t_{0}+s'd_{n}))^{-1}_{km} \\
 && \cdot Y_{il}(t_{0}+s'd_{n}) Y_{im}(t_{0}+s'd_{n})
\lambda_{i}(t_{0}+s'd_{n}) ds'.
\end{eqnarray*}

Since 
\begin{eqnarray*}
\lambda_{i}(t_{0}+s'd_{n}) = 
\sum_{g=0}^{p} \beta_{g}(t_{0}+s'd_{n}) Z_{ig}(t_{0}+s'd_{n}) Y_{i}(t_{0}+s'd_{n}),
\end{eqnarray*}
this becomes
\begin{eqnarray*}
\langle \tilde{v}_{n,j}^{(2)}, \tilde{v}^{(2)}_{n,k} \rangle (s) &=& 
\frac{1}{n} \sum_{i=1}^{n}
\int_{0}^{s} 
\sum_{l,m=0}^{p}
(r^{(2)}(t_{0}+s'd_{n}))^{-1}_{jl}
(r^{(2)}(t_{0}+s'd_{n}))^{-1}_{km} \\
&&\cdot Y_{il}(t_{0}+s'd_{n}) Y_{im}(t_{0}+s'd_{n})
\sum_{g=0}^{p} \beta_{g}(t_{0}+s'd_{n}) Y_{ig}(t_{0}+s'd_{n})ds'.
\end{eqnarray*}
Since, by Assumption \ref{Ass:2},
\begin{eqnarray*}
\sup_{s'\in [0,s]}|\frac{1}{n} \sum_{i=1}^{n}
Y_{il} Y_{im} Y_{ig} (t_0+s'd_n)- r^{(3)}_{lmg}(t_0+s'd_n)| &\stackrel{P}{\to}&0,
\end{eqnarray*}
we obtain
\begin{eqnarray*}
\langle \tilde{v}_{n,j}^{(2)}, \tilde{v}^{(2)}_{n,k} \rangle (s) 
&\stackrel{P}{\to}  &\int_{0}^{s}
\sum_{g,l,m=0}^{p}
(r^{(2)}(t_{0}))^{-1}_{jl}
(r^{(2)}(t_{0}))^{-1}_{km}
r^{(3)}_{lmg}(t_{0}) \beta_{g}(t_{0}) ds'.
\end{eqnarray*}

The result also shows that the asymptotic covariance of $ \tilde{v}_{n,j}^{(2)} $ and $ \tilde{v}^{(2)}_{n,k}$ is given by the right hand side of the above expression. Similarly, if we let $\tilde{v}^{(2)}$ denote the process obtained as the limit in distribution of $\tilde{v}_{n}^{(2)}$, using the conditional indepedence of the increments of the martingale difference sequence, we can show that for $s'',s'''>0$,
\begin{eqnarray} \label{eq:covariance-of-tildev_^2}
   &&Cov(\tilde{v}_j^{(2)}(s''),\tilde{v}_k^{(2)}(s'''))= \nonumber\\
   && \int_{0}^{\min(s'',s''')}
\sum_{g,l,m=0}^{p}
(r^{(2)}(t_{0}))^{-1}_{jl}
(r^{(2)}(t_{0}))^{-1}_{km}
r^{(3)}_{lmg}(t_{0}) \beta_{g}(t_{0}) ds'.
\end{eqnarray}

We next verify the Lindeberg condition for $\tilde{v}_{n}^{(2)}$, for establishing the asymptotic normality. We have, with the choice $d_n=n^{-1/3}$, and $ \tilde{v}_{n,j}^{(2)} $ the $j'$th component, defined in $(\ref{eq:v_n-tilde-j})$,
\begin{eqnarray*}
&& \sum_{j=1}^{r_{n}} \mathbb{E} \left[ 
( \tilde{v}_{n,j}^{(2)} (s))^{2}  \cdot \mathbf{1}_{ \{ 
| \tilde{v}_{n,j}^{(2)}(s) | > \varepsilon \}}  \right]  \leq \\
&& \mathbb{E} \left[
\left( \displaystyle\frac{1}{n} \sum\limits_{i=1}^{n} \int_{0}^{s} 
\sum_{l=0}^{p} 
\left( r^{(2)}(t_{0}+s'd_{n}) \right) ^{-1}_{jl}  
Y_{il}(t_{0}+s'd_{n})
\right)^{2}  \lambda_{i}(t_{0}+s'd_{n}) ds'\cdot  \mathbf{1}_{ \{ 
| \tilde{v}_{n,j}^{(2)}(s) | > \varepsilon \}}  \right]  \\
&\leq&\sum\limits_{i=1}^{r_{n}}
 \dfrac{1}{n}
 \sum\limits_{i=1}^{n} \int_{0}^{s} 
\sum\limits_{l=0}^{p} \sum\limits_{k=0}^{p} 
\sup\limits_{s' \in [0, s]} \left| r^{(2)}(t_{0}+s'd_{n}) \right|^{-1}_{jl} \sup\limits_{s' \in [0, s]}\left| r^{(2)}(t_{0}+s'd_{n}) \right|^{-1}_{jk} \\
&&  \sup\limits_{i= 1,...,n, \, s' \in [0, s]} \left| Y_{il}(t_{0}+s'd_{n})\right|    \lambda_{i}(t_{0}+s'd_{n}) ds' \times \mathbb{E} \left(
\mathbf{1}_{ \left\{ 
\left|  \tilde{v}_{n,j}^{(2)}(s) \right| > \varepsilon \right\}} \right)
\end{eqnarray*} 
where the last inequality follows by expanding the square and the triangle inequality.

Under Assumption \ref{Ass:1},  $r^{(2)}$ are continuous functions on the compact $[0,s]$, and thus they are bounded. 

Furthermore, from Chebyshev's inequality we get,
\begin{eqnarray*}
&&P \left(  \left| \frac{d_{n}^{-2}}{\sqrt{n}} \frac{1}{\sqrt{n}} \sum\limits_{i=1}^{n} \int_{0}^{s}  \sum\limits_{l=0}^{p} \left( r^{(2)}(t_{0}+s'd_{n}) \right)^{-1}_{jl}  Y_{il}(t_{0}+s'd_{n}) \, dM_{i}(t_{0}+s'd_{n}) \right| > \varepsilon 
\right)\\
&\leq& 
 \varepsilon^{-2} E \left[  \displaystyle\frac{1}{n} 
\sum\limits_{i=1}^{n} \int_{0}^{s}  \left( \sum\limits_{l=0}^{p} 
\left( r^{(2)}(t_{0}+s'd_{n})\right)^{-1}_{jl}  Y_{il}(t_{0}+s'd_{n}) \right)^{2} \lambda_{i}(t_{0}+s'd_{n}) ds'  \right] .
\end{eqnarray*}

Thus, our proof is accomplished due to the uniform convergence in probability to zero of the $Y_{il}$ function, in Assumption \ref{Ass:2}.

$(iii)$: {\em The term $\tilde{v}_{n}^{(3)}$ is asymptotically negligible}. \\
From the definition of $\tilde{v}_{n}^{(3)}$, we see that
\begin{eqnarray*}
\tilde{v}_{n}^{(3)}(s)&=&d_{n}^{-1} \int_{0}^{s} \big(J(t_{0}+s'd_{n}) - 1\big) \beta(t_{0}+s
'd_{n})\, ds'.
\end{eqnarray*}

Recall that the process $J(t_{0}+s'd_{n})$ is the indicator that $Y^{T}(t_{0}+s'd_{n})$ has full rank. 
Define the set
\begin{eqnarray*}
E_{n} = 
\left\{
\sup_{s' \in [0,s]} 
\Big\|
\tfrac{1}{n} R^{(2)}(t_{0}+s'd_{n}) - r^{(2)}(t_{0}+s'd_{n})
\Big\| < \varepsilon
\right\}.
\end{eqnarray*}
We have established that $r^{(2)}(t)$ is nonsingular, and hence $\tfrac{1}{n}R^{(2)}(t_{0}+s'd_{n})$ is invertible on $E_{n}$, 
and therefore $J(t_{0}+s'd_{n}) = 1$ for all $s' \in [0,s]$.

Thus, with the choise $d_n=n^{-1/3}$,
\begin{eqnarray*}
&&P \left(
\sup_{s' \in [0,s]} 
\left|
{d_n}^{-1} \int_{0}^{s} 
\big(J(t_{0}+s'd_{n}) - 1\big) \beta(t_{0}+s'd_{n})\, ds
\right| > \varepsilon
\right)  \\
&=&P \left(
\sup_{s' \in [0,s]} 
\left|
{d_n}^{-1} \int_{0}^{s} 
\big(J(t_{0}+s'd_{n}) - 1\big) \beta(t_{0}+s'd_{n})\, ds'
\right| > \varepsilon \cap E_{n}
\right) \\
&&+ P \left(
\sup_{s' \in [0,s]} 
\left|
{d_n}^{-1} \int_{0}^{s} 
\big(J(t_{0}+s'd_{n}) - 1\big) \beta(t_{0}+s'd_{n})\, ds'
\right| > \varepsilon \cap E_{n}^{c}
\right) \\
&\le& 
P(E_{n}^{c})\\
 &\to& 0,
\end{eqnarray*}
where the inequality follows since the first term vanishes since $J(t_{0}+s'd_{n}) = 1$ on $E_n$. Consequently,
\begin{eqnarray*}
P \left(
\sup_{s' \in [0,s]} 
|\tilde{v}_{n}^{(3)} | > \varepsilon
\right) &\longrightarrow & 0,
\end{eqnarray*}
as $n\to \infty$, i.e. $\tilde{v}_{n}^{(3)}$ is asymptotically negligible.

\eop

\noindent
The coordinate-wise result now follows by the Cramér-Wold device. Recall that 
\begin{eqnarray*}
 \tilde{v}_{k,n}(s)&=&d_{n}^{-2} \big( v_{k,n}(t_{0} + s d_{n}) - v_{k,n}(t_{0}) \big),
\end{eqnarray*}
is the coordinate-wise rescaled process, for $k=0,\ldots,p$.
\bcor\label{corr:limit-ditribution-rescaled}
Suppose that Assumptions \ref{Ass:1}- \ref{Ass:3} hold. Then, for any $k=0,\ldots,p$,
\begin{eqnarray*}
      \tilde{v}_{k,n}(s)&\stackrel{d}{\to}& \tilde{v}_{k}(s),
\end{eqnarray*}
on $D(-c,c)$, as $n\to\infty$, where $\tilde{v}$ is mean zero Gaussian process with covariance structure 
\begin{eqnarray*}
Cov( \tilde{v}_{k}(s''), \tilde{v}_{k} (s''')&=& \sigma^2_k\,\min(s'',s'''),
\end{eqnarray*}
where
\begin{eqnarray*}
   \sigma^2_{k}&=&\sum_{g,l,m=0}^{p}
(r^{(2)}(t_{0}))^{-1}_{kl}
(r^{(2)}(t_{0}))^{-1}_{km}
r^{(3)}_{lmg}(t_{0}) \beta_{g}(t_{0}).
\end{eqnarray*}
\ecor
\noindent
\prf
The result follows from the previous theorem and the Cramér-Wold device with choice of coefficients $\alpha_k=1$, and $\alpha_j=0$ for $j\neq k$.
\eop
\noindent
The corollary thus establishes Assumption A1 of \cite{anevski.hossjer.2006} for $ \tilde{v}_{k,n}$. 

\bnote
We note that the limit process $\tilde{v}_k$ can be identified with a (two-sided) Brownian motion with covariance structure $Cov( \tilde{v}_{k}(s'), \tilde{v}_{k} (s'')= \sigma^2_k\,\min(s',s'')$, with $\sigma^2_k$ defined above.  \eop
\enote

Next we define $k$'th coordinate of the rescaled deterministic part
\begin{eqnarray*}
g_{k,n}(s) &=& d_n^{-2} \left( \int_{t_0}^{t_0 + s d_n} Y^{-1}(u) Y(u) \, d B(u) - Y^{-1}(t_0) Y(t_0) \, d B(t_0) sd_{n} \right)_k\\
&=&d_n^{-2} \left( \int_{t_0}^{t_0 + s d_n} \, d B_k(u) - \beta_k(t_0) sd_{n} \right)\ .
\end{eqnarray*}
We will first show that $g_{k,n}$ satisfies Assumption A2 of \cite{anevski.hossjer.2006}, i.e.\;that  for every finite $c>0$ there is an $A_k < 0$ such that
\begin{eqnarray}
   \sup_{|s| \leq c} \left| g_{k,n}(s) - A_k s^2 \right| &\to&  0 \label{Ass:A2}
\end{eqnarray}
as $n \to \infty$. But it is elementary so see that this holds if $\beta_k$ is differentiable with $\beta_k'<0$ in a neighbourhood around $t_0$, with $ A_k = \beta'_k(t_0)$.

Next we want to establish Proposition 1 of  \cite{anevski.hossjer.2006}, from which Assumptions A3 and A4 of that paper will follow.

\blem \label{lem:limit-tail-bound}
Suppose that Assumptions \ref{Ass:1}, \ref{Ass:2} and \ref{Ass:3} hold and that  $\beta_k$ is differentiable with $\beta_k'<0$ in a neighbourhood around $t_0$. Then Proposition 1 in \cite{anevski.hossjer.2006} holds for $\tilde{v}_{k,n}$ and $g_{k,n}$, i.e.\;they  satisfy
$\forall \varepsilon, \delta > 0$, $\exists \tau = \tau(\varepsilon, \delta)<\infty$, such that
\begin{eqnarray*}
 \limsup_{n \to \infty} \mathbb{P} \left( \sup_{|s| \geq \tau} \left| \frac{\tilde{v}_{k,n}(s)}{g_{k,n}(s)} \right| > \varepsilon \right) &<& \delta, 
\end{eqnarray*}

\elem
\noindent
\prf
We show the result by first bounding $g_{k,n}$, and then using that bound to prove, via Doob's and Chebyshev's inequalities and the Ito isometry and properties of $\tilde{v}_{k,n}$, the full result.

$(i)$ {\em Bounding $g_{k,n}(s)$}: We have shown that $g_{k,n}$ satisfies Assumption A2 in \cite{anevski.hossjer.2006}, i.e.\;that $(\ref{Ass:A2})$ holds. Then in  particular, $ \forall \tau > 0$, $0 < \varepsilon < \displaystyle\frac{1}{2} |A_k| \tau^{2} $ and for $s = \pm \tau$, we get
\begin{eqnarray*}
 g_{k,n}(\pm \tau) &\leq& A_k \tau^2 + \varepsilon.
  \end{eqnarray*}
Since $g_{k,n}(0) = 0$ and $g_{k,n}$ is concave, for some finite $n_{0} = n_{0}(\varepsilon) $, we have that 
\begin{eqnarray*}
 g_{k,n}(s) &\leq& \displaystyle\frac{g_{k,n}(\tau)}{\tau} |s| \\
 &\leq& \displaystyle\frac{A_k \tau^2 - \varepsilon}{\tau} |s| \\
 &\leq & \displaystyle\frac{1}{2} A_k \tau^{2-1} |s|,
 \end{eqnarray*}
for all $|s| \geq \tau$ and all $n \geq n_0$. Thus we have established that for all $|s| \geq \tau$ and all $n \geq n_0$, 
\begin{eqnarray}
 g_{n}(s) \geq \displaystyle\frac{1}{2} A_k \tau  |s|. \label{eq:gn_bound}
 \end{eqnarray}

$(ii)$ {\em Bounding $\tilde{v}_{k,n}(s)$:} The proof is similar to the corresponding result for the rescaled process in  \cite{anevski.hossjer.2006} in the cases of rescaled partial sum processes and empirical process, for which one partitioned $\{|s|\geq \tau\}$ into intervals, exhibited bounds of the process at the boundaries of those intervals, and used a modulus of continuity for the processes on the intervals.  We, however, will use the fact that we have martingales to our advantage by using Doob's maximal $L^2$ inequality to bound the maximum over an interval by the values at its endpoints, and then the Ito isometry. 

Thus we partition the tail set $\{|s|\geq \tau\}$ into dyadic intervals, by
\begin{eqnarray*}
    \{|s|\geq \tau\}&=&\cup_{j=0}^{\infty} B_j,
\end{eqnarray*}
with $B_j=\{s:2^j \tau \leq s \leq 2^{j+1} \tau\}$, for $j=0,1,2,\ldots.$  Then for $|s|\geq \tau$, and using the bound $(\ref{eq:gn_bound})$, we get
\begin{eqnarray*}
   \frac{|\tilde{v}_{k,n}(s)|}{|g_{k,n}(s)|} &\leq & \frac{2}{A_k} \frac{ |\tilde{v}_{k,n}(s)|}{\tau  |s|}.
\end{eqnarray*}
Therefore
\begin{eqnarray}
   \mathbb{P}\left(\sup_{|s|\geq \tau}  \frac{|\tilde{v}_{k,n}(s)|}{|g_{k,n}(s)|} > \epsilon \right) &\leq &\mathbb{P} \left( \sup_{|s|\geq \tau}  \frac{|\tilde{v}_{k,n}(s)|}{|s|}  > \frac{A_k}{2} \epsilon \tau \right) \nonumber \\
   &\leq & \mathbb{P} \left( \cup_{j=0}^{\infty} \{ \sup_{s \in B_j}   \frac{|\tilde{v}_{k,n}(s)|}{|s|}  > \frac{A_k}{2} \epsilon \tau \}\right)  \nonumber \\
   &\leq &   \sum_{j=0}^{\infty}   \mathbb{P} \left(  \sup_{s \in B_j} \frac{|\tilde{v}_{k,n}(s)|}{|s|}  > \frac{A_k}{2} \epsilon \tau \right) \nonumber \\
   &\leq &   \sum_{j=0}^{\infty}    \mathbb{P} \left( \sup_{s \in B_j} |\tilde{v}_{k,n}(s)|  >  \epsilon \frac{A_k}{2} \tau^2 2^j \right), \label{eq:bounding-partition}
\end{eqnarray}
where the last inequality follows since on $B_j$ we have $|s|> 2^j \tau$.

We now bound the individual terms in the above sum,  by
\begin{eqnarray}
       \mathbb{P} \left( \sup_{s \in B_j} |\tilde{v}_{k,n}(s)|  >  \epsilon A_k \tau^2 2^{j-1} \right) &\leq & \frac{\mathbb{E} \left( \sup_{s \in B_j} \tilde{v}^2_{k,n}(s) \right)}{( \epsilon A_k \tau^2 2^{j-1})^2} \nonumber \\
       &\leq & \frac{ 4 \mathbb{E} \left( \tilde{v}^2_{k,n}(2^{j+1}\tau) \right)}{( \epsilon A_k \tau^2 2^{j-1})^2}, \label{eq:bounding-individual}
\end{eqnarray}
where the first inequality follows by Chebyshev's inequality and the second by Doob's maximal $L^2$ inequality. 
By the Ito isometry
\begin{eqnarray*}
 \mathbb{E} \left( \tilde{v}^{2}_{k,n}(2^{j+1}\tau) \right) = d_n^{-4} (\int_{t_0}^{t_0 + 2^{j+1}\tau d_n} \| Y^{-}(u) \|^2 \, d\langle M \rangle(u))_k,
 \end{eqnarray*}
which is bounded, by $C<\infty$ say, since $Y^{-}$ is bounded in probability and  $M$ is a square-integrable martingale.

Thus, from $(\ref{eq:bounding-partition})$ and $(\ref{eq:bounding-individual})$, we get
\begin{eqnarray*}
   \mathbb{P}\left(\sup_{|s|\geq \tau}  \frac{|\tilde{v}_{k,n}(s)|}{|g_{k,n}(s)|} > \epsilon \right) &\leq & \frac{64\, C}{(\epsilon A_k \tau^2)^2}\sum_{j=0}^{\infty} 2^{-j}\\
   &<&\delta,
   \end{eqnarray*}
   where the last inequality follows by choosing $\tau=\tau(\epsilon,\delta)$ large enough for fixed $\epsilon,\delta>0$, and $n\geq n_0$.
\eop

Finally we establish the tail behaviour of the components of the limit process $\tilde{v}$, i.e.\,we prove that  $\tilde{v}_k$ satisfies Assumption A5 in \cite{anevski.hossjer.2006}.
\blem\label{lem:limitprocess-bound}
Suppose that Assumptions \ref{Ass:1}, \ref{Ass:2} and \ref{Ass:3} hold and that  $\beta_k$ is differentiable with $\beta_k'<0$ in a neighbourhood around $t_0$, for a fixed $k=0,1,\ldots,p$.  Then the component $\tilde{v}_{k}$ of the limit process $\tilde{v}$ satisfies Assumption A5 in \cite{anevski.hossjer.2006}, i.e.\,for every $\epsilon,\delta>0$, $\exists \tau = \tau(\varepsilon, \delta)<\infty$, such that
\begin{eqnarray*}
 \mathbb{P}\left(\sup_{|s|\geq \tau}  \frac{|\tilde{v}_{k}(s)|}{s^2} > \epsilon \right) &<&\delta.
   \end{eqnarray*}
\elem
\prf
The proof is a straight-forward adaptation of the methods in the proof of Lemma \ref{lem:limit-tail-bound}, with the use of Doob's and Chebyshev's inequalities and the Ito isometry.
\eop

We are next able to state a limit distribution result for the order restricted estimator $\tilde{B}_k$, defined in $(\ref{eq:cumulative-isotonic-regression})$, of the cumulative function $B_k$.

\bth \label{thm:limit-distyribution-cumulative} Suppose that Assumptions \ref{Ass:1}, \ref{Ass:2} and \ref{Ass:3} hold and that  $\beta_k$ is differentiable with $\beta_k'<0$ in a neighbourhood around $t_0$.  Then
\begin{eqnarray*}
      n^{2/3} c(t_0)(\tilde{B}_k(t_0)-B_k(t_0))&\stackrel{d}{\to}& S(-s^2+w(s))(0),
  \end{eqnarray*}
  as $n\to \infty$, where 
  \begin{eqnarray*}
         c(t_0)&=& 2^{-1/3}|\beta_k'(t_0)|^{1/3} (\sigma_k^2)^{-2/3},
  \end{eqnarray*}
  and $w$ is a standard two-sided Brownian motion.
  \eth
  
  \prf
Since we have established that Assumption A1-A5 in \cite{anevski.hossjer.2006} hold, we have that 
    \begin{eqnarray*}
        n^{2/3} [\tilde{B}_k(t_0)-B_k(t_0)]&\stackrel{d}{\to}& S(A_k s^2+\tilde{v}_k(s))(0),
  \end{eqnarray*}
as $n\to \infty$,  as a consequence of Theorem 1 in \cite{anevski.hossjer.2006}. 

Furthermore,  since $\tilde{v}_k$ is a two-sided Brownian motion, with $\tilde{v}_k(0)=0$, and with covariance $Cov(\tilde{v}_k(s),\tilde{v}_k(s'))=\sigma_k^2 \min(s,s')$, we have by the self similarity properties of Brownian motion that $\tilde{v}_k(s)\stackrel{d}{=} (\sigma_k^2)^{1/2}\, w(s)$, with $w$ a standard (two-sided) Brownian motion. In fact, we can simplify the expression for limit distribution further, by the change of variable $s=\gamma u$, to obtain 
    \begin{eqnarray*}
A_k s^2+ \tilde{v}_k(s)&=&A_k \gamma^2 u^2+\tilde{v}_k(\gamma u)\\
& \stackrel{d}{=} &A_k \gamma^2 u^2+\gamma^{1/2} (\sigma_k^2)^{1/2}\,w(u)\\
 &=&(\sigma_k^2)^{2/3}A_k^{-1/3} [-u^2+w(u)],
  \end{eqnarray*}
where the second equality follows by the self similarity of Brownian motion, and the third by choosing $\gamma$ so that $-A_k \gamma^2 =\gamma^{1/2} (\sigma_k^2)^{1/2}$, i.e.\,with $\gamma=-(\sigma_k^2)^{1/3}(-A_k)^{-2/3}$. 

Finally, we use that $S(c g(u))=c S(g(u))$ for any function $g$ and any constant $c>0$, by properties of the least concave majorant $S$, cf.\,e.g.\,Lemma A1 in \cite{anevski.hossjer.2006} (noting the typo in formula (74) in \cite{anevski.hossjer.2006}; the constant $a$ must be positive), to establish that
    \begin{eqnarray*}
S(A_k s^2+ \tilde{v}_k(s))&\stackrel{d}{=}& (\sigma_k^2)^{2/3}(-A_k)^{-1/3} S(-s^2+w(s)).
  \end{eqnarray*}
  Finally, noting that $-A_k=|\beta_k'(t_0)|/2$, proves the formula for $c_1(t_0)$, and ends the proof of  the theorem.

  \eop
In order to state the final limit distribution, for the solution $\tilde{\beta}_k$, we need to study the limit process $y(s)=-s^2 +w(s)$, and show that it satisfies the assumptions of Proposition 2 in \cite{anevski.hossjer.2006}, with the appropriate analog statements for the least concave majorant, and thus that Assumption A6 in  \cite{anevski.hossjer.2006} for $y(s)$ holds. However, this has in fact been already established for the process $y(s)$ in \cite{anevski.hossjer.2006}. Thus we have the following theorem.

\bth\label{th:limit-distribution-grenander}
Suppose that Assumptions \ref{Ass:1}, \ref{Ass:2} and \ref{Ass:3} hold and that  $\beta_k$ is differentiable with $\beta_k'<0$ in a neighbourhood around $t_0$.  Then
\begin{eqnarray*}
      n^{1/3} c(t_0)(\tilde{\beta}_k(t_0)-\beta_k(t_0))&\stackrel{d}{\to}& S(-s^2+w(s))'(0),
  \end{eqnarray*}
  as $n\to \infty$, where 
  \begin{eqnarray*}
         c(t_0)&=& 2^{-1/3}|\beta_k'(t_0)|^{1/3} (\sigma_k^2)^{-4/3},
  \end{eqnarray*}
  and $w$ is a standard two-sided Brownian motion.
  \eth

\prf
Since we have established that Assumptions A1-A6 in  \cite{anevski.hossjer.2006} hold, from Theorem 2 in  \cite{anevski.hossjer.2006}  it follows that 
  \begin{eqnarray*}
        n^{1/3} [\tilde{\beta}_k(t_0)-\beta_k(t_0)]&\stackrel{d}{\to}& S(A_k s^2+\tilde{v}_k(s))'(0),
  \end{eqnarray*}
as $n\to \infty$. Rescaling and use of self similarity for the Brownian motion as in the proof of Theorem \ref{thm:limit-distyribution-cumulative}, shows the statement of the theorem. 
\eop

\bnote
The limit distribution $S(-s^2+w(s))'(0)$ is a version of the Chernoff distribution $\argmax_{s\in {\mathbb R}}(-s^2+w(s))$, that arises in many cases of nonparametric order restricted inference.
\enote

The assumption that $\beta_k$ are decreasing models, as already noted, diseases with early onset covariate effects. Imposing the assumption that $\beta_k$ is {\em increasing} models instead diseases with {\em late-onset} covariate effects. We suggest a corresponding nonparametric order restricted estimator of $\beta_k$ as
\begin{eqnarray*}\label{eq:def-isotonic-regression-2}
     \bar{\beta}^{(k)}&=&\argmin_{\gamma \in {\cal R}_{L_k}} \sum_{i=1}^{L_k} ( \hat{\beta}_i^{(k)}- \gamma_i)^2,
\end{eqnarray*}
where now, for any integer number $L$, let ${\cal R}_L=\{\gamma \in {\mathbb R}^L: \gamma_1\leq \ldots \leq \gamma_L\}$ is the set of non-decreasing vectors. It is straightforward to show that $\bar{\beta}^{(k)}$ is given by 
\begin{eqnarray*}\label{eq:algorithm-isotonic-regression-2}
       \bar{\beta}_k(t)&=&\frac{d}{dt}(\tilde{S}(\hat{B}_k(t)))
\end{eqnarray*}
where $\tilde{S}$ is the {\em greatest convex minorant} map and $d/dt$ denotes the left hand derivative.  Also, the next result is immediate.

\bth\label{th:limit-distribution-grenander-2}
Suppose that Assumptions \ref{Ass:1}, \ref{Ass:2} and \ref{Ass:3} hold and that  $\beta_k$ is differentiable with $\beta_k'>0$ in a neighbourhood around $t_0$.  Then
\begin{eqnarray*}
      n^{1/3} c(t_0)(\bar{\beta}_k(t_0)-\beta_k(t_0))&\stackrel{d}{\to}& \tilde{S}(s^2+w(s))'(0),
  \end{eqnarray*}
  as $n\to \infty$, where 
  \begin{eqnarray*}
         c(t_0)&=& 2^{-1/3}|\beta_k'(t_0)|^{1/3} (\sigma_k^2)^{-4/3},
  \end{eqnarray*}
  and $w$ is a standard two-sided Brownian motion.
  \eth

\section{Asymptotic minimax lower bound}\label{sec:asymptotic-lower-bound}

In this section we use the van Trees inequality \cite{vantrees.1968} (cf.\,also \cite{gill.levit.applications.1995}) to establish a  minimax lower bound for the inference problem that we treat, adapting the approach used in \cite{anevski.2003.estimating} for the analogous problem of estimating the derivative of a convex density.

We first state a result which has independent interest, and which has apparently not been derived previously, namely an asymptotic minimax lower bound for the estimation of a monotone hasard function $h$.

Let $\mathcal F := \mathcal F(\epsilon,C_1,C_2,t_0,A)$ denote the class of decreasing hazard functions $h$, strictly decreasing on a fixed neighbourhood $A$ of $t_0$, i.e.\,with $h'(t)<0$ for $t\in A$, and satisfying, for fixed $\epsilon>0$ and constants $0<C_1<C_2<\infty$, the boundedness condition
\begin{eqnarray*}
 C_1< \inf_{|t-t_0|<\epsilon}h(t) &\leq& \sup_{|t-t_0|<\epsilon}h(t)<C_2.
\end{eqnarray*}
Define the metric $\rho(f,g)=2\int |f(t)-g(t)|\,dt$ on $\mathcal F$.

\begin{theorem}
\label{thm:lower-bound-0}
For every open set $U$ in $\mathcal F$ there exists a constant $0<c(U)< \infty $ such that
\begin{eqnarray*}
\liminf\limits_{n \to \infty} \inf_{T_n} \sup_{h \in U}n^{2/3} E_{h} \big[T_n(t_0)-h(t_0)\big]^2 &\geq& {c(U)}, 
\end{eqnarray*}
where the infimum is taken over all functions of the data.
\end{theorem}

The above theorem will be a special case of the asymptotic minimax lower bound for an individual component in an Aalen additive hasard model, which we state and prove next.

Thus suppose that $\beta_k$ is a individual component in the model $(\ref{eq:Aalen.model})$, and suppose that $\beta_k \in {\cal F}$, with the proviso that possibly $\beta_k$ may be a not full hasard function. For any fixed $\beta_k$, and any fixed vector of covariates $(z_1,\ldots, z_p)$, let ${\cal L}_{\beta_k}$ be the set of distributions that arise from the Aalen model $(\ref{eq:Aalen.model})$. Thus 
\begin{eqnarray*}
      {\cal L}_{\beta_k} &=&\{F_{\beta_k,\beta_{-k},z}(t,): z \in {\mathbb R}^{p+1}, \beta_{-k} \}.
\end{eqnarray*}
where 
\begin{eqnarray*}
F_{\beta_k,\beta_{-k},z}(t)&=&1-\exp{\Big(-\int_{0}^t \big( \beta_k(u) z_k +\sum_{j \neq k} \beta_j(u)z_j \big)\,du \Big)},
\end{eqnarray*}
where $z=(1, z_1,\ldots,z_p)$ and where $\beta_{-k}$ denotes the function vector $(\beta_0,\ldots,\beta_p)$ excluding the $k$'th coordinate $\beta_k$.
\bth \label{thm:minimax-bound-beta}
For every open set $U$ in $\mathcal F$ there exists a constant $0<c(U)< \infty $ such that
\begin{eqnarray*}
\liminf\limits_{n \to \infty} \inf_{T_n} \sup_{\beta_k \in U} \sup_{F\in {\cal L}_{\beta_k}}  n^{2/3} E_{F} \big[T_n(t_0)-\beta_k(t_0)\big]^2 &\geq& {c(U)}, 
\end{eqnarray*}
where the infimum is taken over all functions of the data.
\eth
\prf
The idea is to for every open neighbourhood $U$ exhibit a parametric subclass inside $U$ and then, using a Bayesian reasoning and a prior on the parameter space, use the van Trees inequality, to obtain a lower bound over the nonparametric class of function $\beta_k$. We need to accommodate that $\beta_k$ does not uniquely determine the distribution. This is accomplished by fixing a vector of functions $\beta_{-k}$ and vector of covariates $z$, which will uniquely determine the distribution, and by exhibiting a lower bound for $\beta_k$ for those fixed $\beta_{-k},z$. Next taking the supremum over all distributions in ${\cal L}_{\beta_k}$ will only increase the maximum risk, which will give our full result.

Since we will take a supremum over all distributions in ${\cal L}_{\beta_k}$, and these are determined by the vector of functions $\beta_{-k}$ and the covariates $(z_1,\ldots,z_p)$, we can choose any $\beta_{-k}$ and $(z_1,\ldots,z_p)$ when working inside the supremum over $\beta_k$. We will make the choice that simplifies our calculations the most; thus we choose $\beta_j \equiv 0$, for all $j\neq k$ and $z_j=1$ for all $j$. Then our full hasard in the Aalen model
 is
 \begin{eqnarray*}
       h(t)&=&\beta_k(t),
 \end{eqnarray*}
 and we have reduced our problem to exhibiting a lower bound for a monotone hasard function. 
 
 Thus, and with the use of $h$ in place of $\beta_k$ in order to simplify the notation, our object of study is 
 \begin{eqnarray*}
 \sup\limits_{h \in U} E_{F}(T_n(t_0)-h(t_0))^2,
 \end{eqnarray*}
where $T_n$ is an arbitrary estimator of $h$ at $t_0$. We construct, for each $h\in U$, a parametric subfamily of hazard functions in $U$ through $h$, and apply the van Trees inequality to this subfamily.

First, let $k$ be a fixed real valued $C^1$ function with support $[-1,1]$, satisfying the normalization
\begin{eqnarray*}
\int_{-1}^{1} k(u)\,du = 0, \qquad k(0)=1. \label{eq:kzero}
\end{eqnarray*}
Since $k'$ is continuous with compact support there is a $\tilde{K}<\infty$ such that $|k'|\leq \tilde K$. Define the parametrised family of hazard functions
\begin{eqnarray*}
h_\theta(t) &=& h(t) + \theta\, k\Big(\frac{t-t_0}{b_n}\Big),
\end{eqnarray*}
where $b_n\downarrow 0$ is a bandwidth, to be specified later. Since $h$ and $k$ are $C^1$, so is $h_\theta$.

Next we will choose a parameter space $\Theta$ such that  $h_\theta'$ is negative  on $A$ for all $\theta \in \Theta$, and thus ensuring $h_\theta \in {\cal F}$. Since $h'(t)<0$ on the compact set $\bar A$,
\begin{eqnarray*}
m := \inf_{t\in \bar A}[-h'(t)] &>& 0 
\end{eqnarray*}
which implies $h'(t)\leq -m<0$, for all $ t\in \bar A$, and so
\begin{eqnarray*}
h_\theta'(t) &=& h'(t)+\frac{\theta}{b_n} k'\Big(\frac{t-t_0}{b_n}\Big) \\
 &\leq&  -m+\Big|\frac{\theta}{b_n}\Big|\tilde K.
\end{eqnarray*}
Therefore, if $|\theta/b_n|\tilde K<m$ then $h_\theta'$ will be negative on $A$. Thus, if we choose the parameter space as 
\begin{eqnarray*}
\Theta &=& (-\delta b_n,\ \delta b_n).
\end{eqnarray*}
with $\delta=m/\tilde K$, then $h_\theta \in {\cal F}$ for every $\theta \in \Theta$.

We next check that the parametric family $\{h_\theta:\theta \in \Theta\}$ lies in  $U$ if $\delta$ is small enough. We have that 
\begin{eqnarray*}
\rho(h,h_\theta) &=& 2\int|h_\theta(t)-h(t)|\,dt \\
 &\leq& \ 2|\theta| \int \Big|k\Big(\frac{t-t_0}{b_n}\Big)\Big|\,dt \\
&=& 2|\theta| \, b_n \int |k(u)|\,du \\
 &\leq& \ 2\delta b_n^2 \int |k(u)|\,du,
\end{eqnarray*}
where the third equality follows by a change of variables $u=(t-t_0)/b_n$ and the fourth inequality by $|\theta|<\delta b_n$. Therefore $\rho(h,h_\theta)$ can be made arbitrarily small by choosing $\delta$ small enough, which implies that $\{h_\theta:\theta \in \Theta\} \subseteq U$ if $\delta$ is small enough.

We now apply the van Trees inequality. Let $T_n$ be an arbitrary estimator of $h(t_0)$, and  let $\lambda$ be any prior density on $\Theta$. Then, for $\delta$ sufficiently small,
\begin{eqnarray*}
\sup_{h \in U} E_h \big(T_n(t_0)-h(t_0)\big)^2 &\geq& \sup_{\theta \in \Theta} E_\theta \big(T_n(t_0)-h_\theta(t_0)\big)^2 \\
&\geq& \int_\Theta E_\theta\big(T_n(t_0)-h_\theta(t_0)\big)^2 \lambda(\theta)\,d\theta,
\end{eqnarray*}
Define $\psi(\theta)= h_\theta(t_0)  = h(t_0)+\theta$. Then by the van Trees inequality \cite{vantrees.1968},
\begin{eqnarray}
\int_\Theta E_\theta\big(T_n(t_0)-h_\theta(t_0)\big)^2 \lambda(\theta)\,d\theta &\geq& \frac{\big(\int \psi'(\theta)\lambda(\theta)\,d\theta\big)^2}{n\int I(\theta)\lambda(\theta)\,d\theta+\tilde I(\lambda)}. \label{eq:vantrees}
\end{eqnarray}
where $I(\theta)$ is the Fisher information under $\theta$, and 
\begin{eqnarray*}
\tilde I(\lambda) &=& \int  \frac{\lambda'(\theta)^2}{\lambda(\theta)}\,d\theta. 
\end{eqnarray*}
We will next bound the three terms on the right hand side of \eqref{eq:vantrees}.

$(i)$ {\em The Fisher information} $I(\theta)$.  The density of $T$ under $h_\theta$ is
\begin{eqnarray*}
f_\theta(t) &=& h_\theta(t)\exp(-\int_0^t h_\theta(s)\,ds),
\end{eqnarray*}
and therefore the score function is
\begin{eqnarray}
\frac{\partial \log f_\theta(t)}{\partial \theta} &=& \frac{k\big(\frac{t-t_0}{b_n}\big)}{h_\theta(t)}-\int_0^t k\Big(\frac{s-t_0}{b_n}\Big)\,ds. \label{eq:score-function}
\end{eqnarray}
Now for $t\notin(t_0-b_n,t_0+b_n)$, we obtain
\begin{eqnarray*}
\int_0^t k\Big(\frac{s-t_0}{b_n}\Big)\,ds &=& b_n \int_{-1}^{1}k(u)\,du \\
& =& \ 0,
\end{eqnarray*}
while for $t\in(t_0-b_n,t_0+b_n)$, we have
\begin{eqnarray*}
\Big|\int_0^t k\Big(\frac{s-t_0}{b_n}\Big)\,ds\Big| &\leq& b_n \int_{-1}^{1}|k(u)|\,du \\ 
&= & O(b_n).
\end{eqnarray*}
This implies that
\begin{eqnarray}
E_\theta\big[\left( \int_0^t k\!\left(\frac{s-t_0}{b_{n}}\right)ds \right)^2 \big] & = & O(b_{n}^2), \label{eq:fisher-bound-1}
\end{eqnarray}
and
\begin{eqnarray}
 E_\theta\big[\frac{k\!\left(\frac{t-t_0}{b_{n}}\right)}{h_\theta(t)}\int_0^t k\!\left(\frac{s-t_0}{b_{n}}\right)ds\big]  &= & O(b_{n}^2), \label{eq:fisher-bound-2}
\end{eqnarray}

Therefore, $(\ref{eq:score-function}), (\ref{eq:fisher-bound-1}), (\ref{eq:fisher-bound-2})$, imply
\begin{eqnarray}
I(\theta) &=& E_\theta\Big[\Big(\frac{\partial \log f_\theta(T)}{\partial \theta}\Big)^2\Big] \nonumber \\
& =&  \int \frac{k^2\big(\frac{t-t_0}{b_n}\big)}{h_\theta(t)^2} f_\theta(t)\,dt \ + \ O(b_n^2). \label{eq:fisher-bound-3}
\end{eqnarray}
We note that for $n$ large enough that $b_n<\epsilon$, we get $(t_0-\epsilon,t_0+\epsilon)\supset(t_0-b_n,t_0+b_n)$ and thus
\begin{eqnarray}
0 < C_1  <  \inf_{|t-t_0|<b_n}h_\theta(t)  <  \sup_{|t-t_0|<b_n}h_\theta(t)  <  C_2. \label{eq:bound-on-h_theta}
\end{eqnarray}
Therefore
\begin{eqnarray}
\int \frac{k^2\big(\frac{t-t_0}{b_n}\big)}{h_\theta(t)^2}f_\theta(t)\,dt &=& \int \frac{k^2\big(\frac{t-t_0}{b_n}\big)}{h_\theta(t)}\exp\Big(-\int_0^t h_\theta(s)\,ds\Big)\,dt \nonumber  \\
&\leq& C_1^{-1}\int k^2\Big(\frac{t-t_0}{b_n}\Big)\,dt \nonumber \\
&=& \ b_n C_1^{-1}\int k^2(u)\,du \nonumber\\
& =:& \ b_n C_3. \label{eq:bound-leading}
\end{eqnarray}
where the second inequality follows by the lower bound in $(\ref{eq:bound-on-h_theta})$ and since the support of $k$ is in $[-1,1]$, and since $h_\theta(t)\geq 0$ we may use the trivial bound $\exp(-\int_0^t h_\theta(s)\,ds)\leq 1$, and the third equality follows by the change of variables $u=(t-t_0)/b_n$.

Combining $(\ref{eq:fisher-bound-3})$ and $(\ref{eq:bound-leading})$, we get  
\begin{eqnarray*}
I(\theta)& \leq & b_{n}\,C_3 + O(b_{n}^2) \\
 &=&  b_{n}\,C_3(1+o(1)),
\end{eqnarray*}
as $ b_n\to 0$,  and hence, for $n$ i.i.d.\,observations,
\begin{eqnarray*}
 I_n(\theta) &=& nI(\theta) \\
  &\leq& nb_nC_3(1+o(1)).
\end{eqnarray*}
as $b_n\to 0$.

$(ii)$ {\em The function $\tilde I(\lambda)$}. We note that if we let $\lambda_0$ be an arbitrary density on $(-1,1)$, and define the prior density $\lambda(\theta)=\frac{1}{a}\lambda_0(\theta/a)$ on a parameter space $\Theta=(-a,a)$, for some $a>0$, then
\begin{eqnarray*}
\tilde I(\lambda) &=& \int_{-a}^a \frac{\lambda'(\theta)^2}{\lambda(\theta)}\,d\theta \\ 
& = & \int_{-a}^a \frac{\big(\frac{1}{a^2}\lambda_0'(\theta/a)\big)^2}{\frac{1}{a}\lambda_0(\theta/a)}\,d\theta \\
& =&  \frac{1}{a^2}\int_{-1}^1 \frac{\lambda_0'(u)^2}{\lambda_0(u)}\,du \\
 &=:& \ \frac{I_0}{a^2}.
\end{eqnarray*}
But since we have already chosen $\Theta=(-\delta b_n,\delta b_n)$, and thus  $a=\delta b_n$, we have in fact
\begin{eqnarray*}
\tilde I(\lambda)& = & \dfrac{I_0}{(\delta^2 b_n^2)}.
\end{eqnarray*}

$(iii)$ {\em The term} $ \int \psi'(\theta)\lambda(\theta)\,d\theta$. Since $\psi(\theta)=h(t_0)+\theta$, we have $\psi'(\theta)=1$, so that 
\begin{eqnarray*}
\left(\int \psi'(\theta)\lambda(\theta)\,d\theta\right)^2 & = & \left(\int \lambda(\theta)\,d\theta\right)^2 \\
& = & 1.
\end{eqnarray*}
Substituting the three bounds obtained in $(i),(ii),(iii)$ into \eqref{eq:vantrees}, we get
\begin{eqnarray*}
\sup_{h \in U} E_h\big(T_n(t_0)-h(t_0)\big)^2 &\geq& \frac{1}{n b_n C_3(1+o(1))+\dfrac{I_0}{\delta^2 b_n^2}}.
\end{eqnarray*}
The two terms in the denominator are of the same order when $nb_n \asymp b_n^{-2}$, i.e.\,when $b_n\asymp n^{-1/3}$. With this bandwidth choice
\begin{eqnarray*}
nb_nC_3 &=& C_3\,n^{2/3}, \qquad \frac{I_0}{\delta^2 b_n^2} \ = \ \frac{I_0}{\delta^2}\,n^{2/3},
\end{eqnarray*}
so that the denominator is $O(n^{2/3})$. Thus, and now we revert to the notation $\beta_k=h$ and $F$ for the corresponding distribution,
\begin{eqnarray*}
n^{2/3}\sup_{\beta_k \in U} E_{F}\big(T_n(t_0)-\beta_k(t_0)\big)^2 &\geq& c(U),
\end{eqnarray*}
for the specific $F$ arising with the choice $\beta_j\equiv 0$ for all $j\neq k$ and $z_j=1$ for all $j$, for $n$ large enough. We note that the constant $c(U)>0$ depends on $C_1,C_3,I_0,\delta$ but not on $T_n$ or on $n$. Since taking a supremum of the left hand side over all $F \in  {\cal L}_{\beta_k}$, inside the supremum over $\beta_k$'s, will only increase the left hand side, we get the result
\begin{eqnarray*}
 \sup_{\beta_k \in U} \sup_{F\in {\cal L}_{\beta_k}}  n^{2/3} E_{F} \big[T_n(t_0)-\beta_k(t_0)\big]^2 &\geq& {c(U)},
\end{eqnarray*}
for $n$ large enough. Finally, noting that the right hand side does not depend on the estimator $T_n$, the proof is complete.
\eop

\section{Discussion}\label{sec:discussion}
In this paper we have derived limit distributions for the {\em coordinate wise} least squares projection of a naive estimator on the space of decreasing functions. The results are derived using a general approach presented in \cite{anevski.hossjer.2006}, and the main work in this paper has been to establish the necessary conditions required in \cite{anevski.hossjer.2006} for the conclusions of that paper to hold. That in fact gives us our two main results Theorems \ref{thm:limit-distyribution-cumulative} and \ref{th:limit-distribution-grenander}. The conditions under which we are able to establish these results are the conditions required in 
 \cite{andersen.borgan.gill.keiding.1993} for the derivation of limit distributions of the starting estimator $\hat{B}$; thus we do not need to demand more than is demanded in \cite{andersen.borgan.gill.keiding.1993}.

 One of main vehicles for this is our Theorem \ref{thm:limit-distribution-rescaled-proc}, which derives the limit distribution for the rescaled process $\tilde{v}_n$. We note that the result in Theorem \ref{thm:limit-distribution-rescaled-proc} is in fact stronger than necessary for our need, and that we only need its consequence Corollary \ref{corr:limit-ditribution-rescaled}.  
 
Furthermore in Theorem \ref{thm:minimax-bound-beta} we provide an asymptotic minimax lower bound for the problem of estimation the individual $\beta_k$'s under order restrictions, which shows that the rate $n^{1/3}$ is in fact the optimal rate for this inference problem.

The proposed estimator is feasible to use for modelling and estimation diseases with early onset, or late onset, covariate effects and thus has a clear clinical interpretation.  The estimator provides estimates of the {\em individual} covariate effects in a {\em multivariate model}, and thus is feasible to use as disease model for several covariates that individually have early onset or late onset. For instance one may consider a disease in which one models the genetic effect as an early onset covariate and an environmental or phenotypic covariate as a late onset covariate, thereby enabling the study of diseases with complex gene-environment causes. 

The estimator is completely data-adaptive and does not need bandwidth choosing. Furthermore, the algorithms needed to calculate the estimator are available in standard R packages, and thus our suggested estimator is readily implementable in a clinical setting.

We would like to finally note that although the Aalen additive hasard model is mainly used in survival analysis in clinical settings, our methods do not depend on this assumption.

\section{Acknowledgments}
The research of DA is partially supported by the Swedish Research Council (SRC). DA gratefully acknowledges the SRC's support. We would like to thank Eduard Belitser for suggesting to study the optimal rate question, leading to the minimax lower bound results in Section \ref{sec:asymptotic-lower-bound}.

\end{document}